\documentclass{amsart}
\usepackage{amsmath,amsthm,verbatim, amscd, epsfig}

\newtheorem{theorem}{Theorem}[section]

\theoremstyle{definition}

\newcommand{\R}{\mathbb{R}}
\newcommand{\C}{\mathbb{C}}
\renewcommand{\H}{\mathbb{H}}
\renewcommand{\P}{\mathbb{P}}
\renewcommand{\Im}{{\rm Im}}
\renewcommand{\Re}{{\rm Re}}
\newcommand{\del}{\partial}
\newcommand{\delbar}{\bar{\del}}
\newcommand{\<}{\langle}
\renewcommand{\>}{\rangle}
\newcommand{\Span}{\mathrm{span}}
\begin{document}

\title{Schwarzian Derivatives and Flows of surfaces}
\author{Francis Burstall}
\address{Department of Mathematics\\
University of Bath\\
England}
\email{feb@math.bath.ac.uk}

\author{Franz Pedit}
\address{Department of Mathematics\\
University of Massachusetts\\
Amherst, MA 01003\\
USA}
\email{franz@gang.umass.edu}
\author{Ulrich Pinkall}
\address{Fachbereich Mathematik, MA 8-5\\
Technische Universit\"at Berlin\\
Strasse des 17. Juni 135\\
10633 Berlin\\
Germany}
\email{pinkall@sfb288.math.tu-berlin.de}

\thanks{Authors supported by SFB288 at TU-Berlin.  First author
partially supported by The European Contract Human Potential Programme,
Research Training Network HPRN-CT-2000-0010.  Second author
additionally supported by NSF grants DMS-90011083 and DMS-9705479}
\maketitle

\section{introduction}\label{sec:introduction}

Over the last decades it has been widely recognized that many
completely integrable PDE's from mathematical physics arise
naturally in geometry. Their integrable character in the geometric
context -- usually associated with the presence of a Lax pair and
a spectral deformation -- is nothing but the flatness condition of
a naturally occurring connection. Examples of interest to
geometers generally come from the integrability equations of
special surfaces in various ambient spaces: among the best
known examples is the sinh-Gordon equation describing constant mean
curvature tori in $3$-space. In fact, most of the classically
studied surfaces, such as isothermic surfaces, surfaces of
constant curvature and Willmore surfaces, give rise to such
completely integrable PDE. A common thread in many of these
examples is the appearance of a harmonic map into some symmetric
space, which is well-known to admit a Lax representation with
spectral parameter. Once the geometric problem is formulated this
way algebro-geometric integration techniques give explicit
parameterizations of the surfaces in question in terms of theta
functions.

One contribution of the geometric view of these integrable PDE is a
much better understanding of the meaning of the hierarchy of flows
associated to these equations. In mathematical physics these
hierarchies are obtained by deformations of the Lax operators
preserving their eigenvalues. This is rather unsatisfactory from the
geometric viewpoint, where one wants to see these flows as geometric
deformations on the surfaces. A good example of the geometric
derivation of the KdV flows can be found in \cite{Pin95} (see also
\cite{Seg91}) and the flows related to the sinh-Gordon equation are
explicitly derived from geometric deformations in \cite{PinSte89}.

Lately two more integrable hierarchies, well known to mathematical
physicists, have made their appearance in surface geometry: these are
the Novikov--Veselov and Davey--Stewartson equations
\cite{KonTai95A,Tai97,Kon00,KonLan00}. They are derived as kernel
dimension preserving deformations of the Dirac operator. Since
surfaces in $3$ and $4$ space can be presented by ``spinors", these
flows in principle give an integrable hierarchy on the space of
conformally parameterized tori in $3$ and $4$ space. The relation to
surface geometry leaves a lot to ask for: what is the geometric
meaning of these flows? Can one derive them from some geometric
principles? What kind of geometry is involved?  Even though existing
literature \cite{Kon00,KonLan00,Tai97} uses Euclidean surface
geometry, it turns out that these flows are really M\"obius
invariant. Thus there should be some M\"obius invariant setting from
which these flows can be derived. However, notwithstanding this, the
essential new insight is that there is a completely integrable
hierarchy on {\em all} conformally immersed $2$-tori, and not just on
special classes given by certain geometric conditions, such as
constant mean curvature tori etc. One clearly expects that these new
flows reduce in some way to the known ones when restricted to the
corresponding special surface classes.

The picture emerging here is that the manifold of conformally
parameterized $2$-tori in $3$-space is stratified by finite
dimensional submanifolds, namely those tori stationary with
respect to some higher flow in the hierarchy. These ``finite type"
tori can then be parameterized by theta functions on some finite
genus Riemann surface. As said, this is just a picture and we feel
there still is a lot more explaining to be done to become a useful
mathematical theory.

This present note aims to give a self contained, low technology
approach to the above mentioned topics. The natural setting of our
discussion is M\"obius invariant surface geometry. Since we do not
claim to have a complete theory, the exposition will at times be
sketchy and leave the Reader, hopefully, with some urge to dwell
further on the issues.

The most popular integrable hierarchy is the KdV hierarchy.
Therefore we begin this paper by deriving this hierarchy as
natural flows on holomorphic maps into the Riemann sphere. The
basic invariant of a holomorphic map $f$, its Schwarzian
derivative $S_z(f)$ with respect to a coordinate $z$, turns out to
be the function satisfying the equations of the KdV hierarchy. The
Schwarzian occurs naturally in a Hill equation
\[
\psi_{zz}+\tfrac{1}{2}S_z(f)\psi=0
\]
for some appropriate homogeneous lift $[\psi]=f$. Notice that
under a change of coordinates the Schwarzian transforms as
\[
S_w(f)dw^2=(S_z(f)-S_z(w))dz^2\,,
\]
which makes it conceptually harder to understand what object the
Schwarzian really is (see, however \cite{BurCal}). This approach to
KdV is well known \cite{Seg91}, but our derivations are closer to the
M\"obius geometry of conformally immersed Riemann surfaces into the
$n$ sphere.

The latter theme gets introduced in the second chapter, where we
derive the invariants, structure and integrability equations of a
conformal surface $f:M\to S^n$ under the M\"obius group. In this
case the basic equation defining the invariants is an
inhomogeneous Hill equation
\[
\psi_{z z}+\tfrac{1}{2}S_z(f)\psi=\kappa\,,
\]
for an appropriate lift $\psi$ of $f$ into the forward light
cone---we view the conformal $n$ sphere as the projectivized light
cone in Lorentz space. The invariant $\kappa$ is essentially the
trace free second fundamental form, or {\em Hopf differential}, of
$f$, a well known M\"obius invariant, and $S_z(f)$ is the {\em
Schwarzian derivative}\footnote{The appearance of Schwarzian
derivatives in conformal surface geometry is not new: they were
introduced in a closely related context by Calderbank \cite{Cal98}
who has also developed, with the first author, an independent and
more invariant formulation of the fundamental theorem of M\"obius
invariant surface theory \cite{BurCal}.} of the surface $f$. It
transforms under coordinate changes just like the classical
Schwarzian above. Of course, if $\kappa=0$ we reduce to the classical
case of holomorphic maps into $S^2$. The Hopf differential and
Schwarzian form a complete set of invariants for surfaces in $3$
space with respect to the M\"obius group. Their integrability
equations, the conformal Gauss and Codazzi equations, are
\begin{gather*}
\tfrac{1}{2}S_z(f)_{\bar{z}}=3\bar{\kappa}\kappa_z+\bar{\kappa}\kappa_z\\
\text{Im}\,(\kappa_{\bar{z}\bar{z}}+\tfrac{1}{2}\overline{S_z(f)}\kappa)=0\,,
\end{gather*}
which gives the fundamental theorem of M\"obius invariant surface
theory.

To demonstrate the effectiveness of our setup, we discuss a number
of natural problems from conformally invariant surface theory:
cyclides of Dupin are characterized by holomorphic Schwarzians
$S_z(f)$; isothermic surfaces by real Hopf differentials $\kappa$;
Willmore surfaces are characterized by
\[
\kappa_{\bar{z}\bar{z}}+\tfrac{1}{2}\overline{S_z(f)}\kappa =0\,,
\]
a stronger condition than the Codazzi equation. From these conditions
we easily see the associated families of isothermic and Willmore
surfaces, crucial to their completely integrable character. Moreover,
unless a surface is isothermic, the Hopf differential alone
determines it uniquely up to M\"obius transformations. We also
discuss results of Thomsen \cite{Tho23} and Richter \cite{Ric97}
characterising minimal, respectively constant mean curvature,
surfaces in $3$-dimensional space forms as isothermic Willmore,
respectively constrained Willmore, surfaces and indicate a
generalisation to isothermic Willmore surfaces in $n$-space.  We
include these sometimes very classical results since they amplify the
``retro" charm of our view point on M\"obius geometry: all
calculations are local and done in a fixed coordinate, just like
$100$ years ago, but due to our non-redundant set of invariants turn
out to be extraordinarily short and efficient.

Having done that much, we finally get to a geometric derivation of
the Novikov--Veselov and Davey--Stewartson equations in the last
chapter. We follow our initial derivation of the KdV hierarchy in
the first chapter to derive the Novikov--Veselov flow for surfaces
in $n$-space. This flow is M\"obius invariant by construction and
preserves the Willmore energy, at least in $3$ space. Moreover, it
also preserves isothermic surfaces in any codimension giving a
vector version of the modified Novikov--Veselov flow. Again, when
the Hopf differential $\kappa=0$, so that we deal with holomorphic
maps into $S^2$, these flows reduce to the KdV flows. Whereas the
Novikov--Veselov flows can be defined in any codimension, the Davey
Stewartson hierarchy needs the additional datum of a complex
structure in the normal bundle, which is the case for surfaces in
$4$ space. This flow is defined in term of double derivatives of
the Hopf differential and thus simpler then the Novikov--Veselov
flow, which is given in triple derivatives of the Hopf
differential. Surfaces in $3$ space are generally moved into $4$
space by the Davey--Stewartson flow, which therefore does not
induce flows on surfaces in $3$ space. Nevertheless, the flow
preserves Willmore surfaces and isothermic surfaces in $4$ space.

Even though we mention the word ``hierarchies" frequently, we do
not really have an explicit scheme to construct all the higher
flows for Novikov--Veselov and Davey--Stewartson. This is partially
due to the inherently non-explicit character of the higher flows
---one needs to solve a $\delbar$-problem for each flow--- and
partially due to our insufficient understanding of invariant
multilinear differential operators on Riemann surfaces. The latter
seems to be an important enough issue to return to. At present
there are two attempts in that direction, one \cite{Boh} using
deformation theory on quaternionic holomorphic line bundles for
surfaces in $S^4=\H\P^1$, and another \cite{BurCal} using Cartan
connections and M\"obius structures for surfaces in $S^n$.

This work began during the conference on ``Integrable systems and
geometry'' in Tokyo in the summer of 2000 and has been completed
during January 2001, when all authors met again at the SFB288 at
TU-Berlin. We take this opportunity to thank the organizers of the
Tokyo conference and the SFB288 for their support.

\section{The classical Schwarzian and the KdV hierarchy}
\subsection{The classical Schwarzian derivative}\label{subsec:classical}
Given a meromorphic function $f$ on a Riemann surface $M$ and a holomorphic
coordinate $z$ the {\em Schwarzian derivative of $f$ with respect to $z$}
is given by the expression
\begin{equation}\label{eq:Schwarzian_f}
S_{z}(f)=(\frac{f''}{f'})'-\tfrac{1}{2}(\frac{f''}{f'})^2\,.
\end{equation}
One can quickly check that the Schwarzian $S_{z}(f)=0$ vanishes if and only if
$f=\tfrac{az+b}{cz+d}$ is fractional linear in $z$. Moreover, two meromorphic functions
$f$ and $g$ are related by a fractional linear transformation
\[
g=\frac{af+b}{cf+d}
\]
if and only if their Schwarzians $S_{z}(f)=S_{z}(g)$ agree. The last property
indicates that the Schwarzian is an invariant for maps into $\C\P^1$.
We will develop this
view point further since it motivates most of
our constructions on conformally immersed surfaces later on.

Let $f:M\to \C\P^1$ be a holomorphic immersion, i.e., $df_p\neq 0$
for all $p\in M$, and fix a holomorphic coordinate $z$ on $M$. A {\em
lift} of $f$ is a holomorphic map $\psi:M\to \C^2\setminus\{0\}\,$ so
that $f(p)$ is the line spanned by $\psi(p)$. In other words, $\psi$
is a holomorphic section of the pull-back via $f$ of the tautological
line bundle over $\C\P^1$.  Our assumption that $f$ is an immersion
is equivalent to $\psi$ and $\psi_z$ being linearly independent at
each point $p\in M$.  Using a fixed complex volume form $\det$ on
$\C^2$ we can find, up to sign, a unique {\em normalized} lift $\psi$
of $f$ by demanding that
\begin{equation}\label{eq:psi_normalization}
\det(\psi,\psi_z)=1\,.
\end{equation}
Indeed, two lifts $\tilde{\psi}=\psi\lambda$ are related by a nowhere
vanishing scale $\lambda:M\to\C$ and
\[
\det(\tilde{\psi},\tilde{\psi}_z)=\lambda^2\det(\psi,\psi_z)\,,
\]
which proves our assertion. Taking the derivative of \eqref{eq:psi_normalization}
we obtain
\[
\det(\psi,\psi_{z z})=0\,,
\]
so that $\psi$ satisfies Hill's equation
\begin{equation}\label{eq:Hill_classical}
\psi_{zz}+\tfrac{c}{2}\psi=0
\end{equation}
for some holomorphic function $c:M\to\C$. From our construction it is clear
that $c$ depends only on $f$ and the coordinate $z$. If the holomorphic
map $f$ into $\C\P^1$ is written in homogeneous coordinates $[f:1]$
then one easily checks that
\begin{equation}\label{eq:c}
c=S_{z}(f)
\end{equation}
is the Schwarzian derivative. The above mentioned properties of the
Schwarzian now follow immediately from \eqref{eq:Hill_classical}.
Moreover, given $c$ we solve Hill's equation for
\[
\psi=\begin{pmatrix} \psi_1 \\ \psi_2\end{pmatrix}\,,
\]
and then $f=[\psi_1:\psi_2]$  gives a map into $\C\P^1$, unique up to the
action of ${\bf Sl}(2,\C)$, with Schwarzian $S_z(f)=c$.
Thus the Schwarzian is an infinitesimal invariant, with respect
to the M\"obius group ${\bf PSl}(2,\C)$, for maps into
$\C\P^1$.

To obtain a better understanding what type of object the Schwarzian is, we differentiate
with respect to another holomorphic coordinate $w$ on $M$. An
elementary calculation using \eqref{eq:Hill_classical} reveals
the transformation formula
\begin{equation}\label{eq:coordinate_change}
S_{w}(f)dw^2=(S_z(f)-S_z(w))dz^2\,.
\end{equation}
Recall that a {\em projective structure} on $M$ is given by a holomorphic atlas
$\mathcal{A}$ whose coordinate transition functions are fractional linear.
Given two projective structures $\mathcal{A}$ and $\tilde{\mathcal{A}}$
the holomorphic quadratic
differential
\[
Q:=S_{z}(\tilde{z})dz^2\,,
\]
where $z$ and $\tilde{z}$ are charts in $\mathcal{A}$ and
$\tilde{\mathcal{A}}\,$, is well-defined by
\eqref{eq:coordinate_change}. Therefore the space of all
projective structures on a Riemann surface is an affine space
modelled on the vector space $H^{0}(K^2)$ of holomorphic quadratic
differentials. If we have fixed a projective structure on $M$, the
transformation rule \eqref{eq:coordinate_change} implies that the
Schwarzian of a function may be regarded as a holomorphic
quadratic differential.

\subsection{The KdV hierarchy}\label{subsec:KdV_hierarchy}
Consider a deformation $f(t):M\to\C\P^1$ of $f=f(0)$ with corresponding
lift $\psi(t)$, which we may assume to be normalized \eqref{eq:psi_normalization}. Then
$\psi$ and $\psi_z$ are a unimodular basis for $\C^2$ for any $z$ and $t$,
and the infinitesimal variation of $\psi$ can be expressed as
\begin{equation}\label{eq:psi_dot}
\psi_t=a\psi+b\psi_z\,.
\end{equation}
Here $a$ and $b$ are functions of $t$ and $z$ which have to satisfy
\begin{equation}\label{eq:a=b_z}
b_z=-2a
\end{equation}
due to \eqref{eq:psi_normalization}. For the infinitesimal
variation of the Schwarzian, we use Hill's equation
\eqref{eq:Hill_classical} and obtain
\[
\psi_{tzz}+\tfrac{c_t}{2}\psi+\tfrac{c}{2}\psi_t=0\,.
\]
Calculating $\psi_{tzz}$ from \eqref{eq:psi_dot} and
\eqref{eq:a=b_z} the $\psi$-part of this relation unravels to
\begin{equation}\label{eq:c_dot}
c_t=b_{zzz}+2b_z c+bc_z\,.
\end{equation}
So far we have been describing general deformations of $f$. To
obtain geometrically meaningful deformations the ``radial
variation" $a$, and therefore the ``tangential variation" $b$, in
\eqref{eq:psi_dot} have to be chosen geometrically. The most naive
choice of variation is {\em no} radial variation at all, i.e.,
$a_1=0$. Then, up to a constant, $b_1=1$ by \eqref{eq:a=b_z} and
from \eqref{eq:c_dot} we get
\[
c_{t_1}=c_z\,,
\]
which is translational flow along $z$. To obtain a more interesting flow
we use the {\em old} variation of $c$ as the {\em new} radial variation,
i.e., we put
\[
a_3=-\tfrac{1}{2}c_{t_1}\,.
\]
Then $b_3=c$ by \eqref{eq:a=b_z} and from \eqref{eq:c_dot} we
obtain
\[
c_{t_3}=c_{zzz}+3cc_z\,,
\]
which we recognize as the KdV equation. Continuing our recursion scheme
\[
a_{2k+1}=-\tfrac{1}{2}c_{t_{2k-1}}\,,
\]
we arrive at an infinite hierarchy of odd order (in z) flows, the KdV hierarchy. For example,
the $5^{th}$-order flow computes to
\[
c_{t_5}=c_{zzzzz}+5c_{zzz}c+10c_{zz}c_z+\tfrac{15}{2}c_z c^2\,.
\]
Notice that at each step we have constants of integration, which amounts to
considering a given order flow modulo the lower order flows.

A map
$f:M\to \C\P^1$ has {\em finite type $n$} if the $(2n-1)^{st}$--order flow
is stationary modulo the lower order flows:
\begin{equation}\label{eq:finite_type}
c_{t_{2n-1}}+\sum_{k=1}^{n-1}\lambda_k c_{t_{2k-1}}=0\,.
\end{equation}
For example, type $1$ maps have $c_z=0$ or  $c=c_0$ constant, and
thus are the loxodromes $f=\exp(i\sqrt{2c_0}z)$. In general,
finite type $n$ maps are given by $\theta$--functions on a genus
$n-1$ Riemann surface, and the flows in the hierarchy become
linear flows on the Jacobian of this Riemann surface. The above
geometric construction of an integrable hierarchy and its
algebraic-geometric integration is folklore by now. Our exposition
follows the spirit of \cite{PinSte89}, where a similar construction was
done for the $\sinh$--Gordon hierarchy.

We conclude this paragraph with a remark on the closely related
mKdV hierarchy. Whereas the KdV flows commute with the action of
M\"obius group
${\bf PSl}(2,\C)$, the mKdV flows arise
by ``symmetry breaking": fixing a point at infinity on $\C\P^1$
we consider $f:M\to \C\subset \C\P^1$ where $\C$ now carries the flat
Euclidean structure. In other words, we have reduced the symmetry group
of our target space from the M\"obius group to the Euclidean group.
The fundamental Euclidean invariant of $f$ is now given by the induced metric
\[
|df|^2=f_z\overline{f_z}|dz|^2\,,\qquad  f_z\overline{f_z}=e^{2v}
\]
where $v:M\to \R$ is a real valued function.
Using the formula \eqref{eq:Schwarzian_f} for $S_z(f)$, it
is easy to check that
\begin{equation}\label{eq:Schwarzian_Euclidean}
S_z(f)=2v_{zz}-2(v_z)^2\,.
\end{equation}
In particular, this combination of derivatives of the real valued function
$v$ is holomorphic.
Putting $p:=2v_z$ the above relation reads
\[
S_z(f)=p_z-\tfrac{1}{2}p^2\,,
\]
which we recognize as the {\em Miura transform} \cite{Mum91}. Therefore
the Miura transform is just the expression of the M\"obius invariant
Schwarzian of $f$ in terms of the Euclidean invariant induced metric
of $f$. It is a well known fact \cite{Mum91} that the Miura transform
interpolates between the mKdV and KdV flows so that $p$ satisfies
the mKdV hierarchy.

\section{Conformally immersed surfaces and Schwarzians}
The constructions of the previous section carry over almost
verbatim to conformally immersed surfaces in $n$-space. As a
result we obtain an elegant description of M\"obius invariant
surface geometry. Later we will indicate how to construct M\"obius
invariant flows on conformally immersed surfaces.

\subsection{The light-cone model}
\label{sec:light-cone-model}

Following Darboux \cite{Dar87}, we linearize the action of the M\"obius
group on the sphere $S^n\subset\R^{n+1,1}$ via the diffeomorphism
\begin{equation}\label{eq:light_cone}
S^n \cong \P(\mathcal{L}): x \leftrightarrow [1:x]
\end{equation}
between the sphere and the projectivized light cone where
$\mathcal{L}\subset \R^{n+1,1}$ is the null cone in
$(n+2)$--dimensional Minkowski space for the quadratic form $\langle
v,v\rangle =-v_0^2+\sum_{k=1}^{n+1}v_k^2$.  The projective action of
the Lorentz group on $\P(\mathcal{L})$ is by conformal
diffeomorphisms giving rise to a double covering of the M\"obius
group by the Lorentz group.

The isomorphism \eqref{eq:light_cone} is a special case of a more
general construction which realizes all space forms as conic sections
of $\mathcal{L}$: for non-zero $v_0\in\R^{n+1,1}$, set
\begin{equation*}
S_{v_0}=\{v\in\mathcal{L}: \<v,v_0\>=-1\}.
\end{equation*}
$S_{v_0}$ inherits a positive definite metric of constant
curvature $-\<v_0,v_0\>$ from $\R^{n+1,1}$ and is a copy of a
sphere, Euclidean space or ball according to the sign of
$-\<v_0,v_0\>$.  For \eqref{eq:light_cone}, $v_0=e_0$.

In this picture, it is easy to describe the $m$-spheres (totally
umbilic submanifolds) in $S^n$ \cite{Roz48}: they are all given by a
decomposition $\R^{n+1,1}=V\oplus V^\perp$ with $V^\perp$ space-like
of dimension $(n-m)$.  The $m$-dimensional sphere is then
$\P(\mathcal{L}\cap V)$.  Viewed as a submanifold $S_{v_0}\cap V$
of the conic section $S_{v_0}$, one readily shows that this
$m$-sphere has mean curvature vector
\begin{equation}
\label{eq:1}
H_v=-v_0^\perp-\<v_0^\perp,v_0^\perp\>v
\end{equation}
at $v$, where $v_0=v_0^T+v_0^\perp$ according to the decomposition
$V\oplus V^\perp$.

\subsection{The Schwarzian derivative of a conformally immersed surface}
Let $f:M\to S^n\subset \R^{n+1}$ be a conformal immersion of a
Riemann surface $M$.  A {\em lift} of $f$ is a map $\psi:M\to
\mathcal{L}_{+}$ into the forward light cone such that the null line
spanned by $\psi(p)$ is $f(p)$.  Perhaps the most naive choice for a
lift is the {\em Euclidean lift}
\begin{equation}\label{eq:lift_Euclidean}
\phi=(1,f)\,,
\end{equation}
but any positive scale of this lift will do. Note that the induced metric of
a scaled lift $\phi\lambda$, with $\lambda:M\to \R_{+}$,
is given by
\[
|d(\phi\lambda)|^2=\lambda^2|df|^2\,,
\]
where the induced metric $|df|^2$ is computed with respect to the
standard round metric on $S^n$ which we also denote by $\langle
\,,\,\rangle $.  Therefore the various lifts of $f$ give rise to the
various metrics in the conformal class of $|df|^2$.  From this we see
that there is a unique {\em normalized} lift $\psi$, with respect to
a given holomorphic coordinate $z$ on $M$, for which the induced
metric
\begin{equation}\label{eq:normalized_lift}
|d\psi|^2=|dz|^2
\end{equation}
is the standard flat metric of the coordinate $z$. Since
\eqref{eq:normalized_lift} is invariant under Lorentz
transformations the normalized lift $\psi$ is M\"obius invariant.

If we pass to another holomorphic coordinate $w$ the normalized lift
with respect to this new coordinate
\begin{equation}\label{eq:lift_change}
\tilde{\psi}=\psi|w_z|
\end{equation}
is obtained by scaling $\psi$ by the length of
the Jacobian of the coordinate transformation.
This follows immediately from the fact that $\psi$ is null
and therefore
\[
|d\tilde{\psi}|^2=|\,d\psi|w_z|+\psi d|w_z|\,|^2=|d\psi|^2|w_z|^2=
|dz|^2|w_z|^2=|dw|^2\,.
\]

A fundamental construction of M\"obius invariant surface geometry
is the {\em mean curvature} or \emph{central sphere congruence}
\cite{Tho23} which assigns to each
point $f(p)$ on the surface the unique 2-sphere $S(p)$ tangent to
$f$ at that point and with the same mean curvature vector
$H_{S(p)}=H_p$ at $f(p)$ as $f$.  While this is phrased in terms of
the Euclidean lift $\phi=(1,f):M\to S_{e_0}$, the construction is
conformally invariant: $S(p)$ corresponds to a decomposition
$\R^{n+1,1}=V(p)\oplus V(p)^\perp$ with
\[
V(p)=\Span\{\phi(p),d\phi_p, e_0^T\}
\]
since $\phi$ is tangent to $V(p)\cap S_{e_0}$ at $p$.  In view of
\eqref{eq:1}, this reads
\begin{equation*}
\Span\{(1,f(p)),(0,df_p),(-1,-H_p)\}=\Span\{(1,f(p)),(0,df_p0),(0,H_p-f(p))\}.
\end{equation*}
However, one readily checks that
$(0,H-f)\;\parallel\;\phi_{z\bar{z}}$ for any holomorphic coordinate
$z$ so that $V$ is a rank $4$ subbundle of $M\times\R^{n+1,1}$ of
signature $(3,1)$ given by
\begin{equation}\label{eq:mean_sphere}
V=\text{span}\{\phi,d\phi,\phi_{z\bar{z}}\}.
\end{equation}
It is clear that $V$ so defined is independent of all choices (of
lift or holomorphic coordinate) and so is a conformal invariant.

To facilitate calculations we choose the unique section
$\hat{\psi}\in \Gamma(V)$ with
\begin{equation}\label{eq:psi_hat}
\langle \hat{\psi},\hat{\psi}\rangle =0\,,\qquad
\langle \psi,\hat{\psi}\rangle =-1\,,\qquad\text{and} \qquad
\langle \hat{\psi},d\psi\rangle =0\,.
\end{equation}
Then, given a holomorphic coordinate $z$, we have the M\"obius
invariant framing
\begin{equation}\label{eq:frame_V}
\psi \,,\;\psi_z\,,\; \psi_{\bar{z}}\,,\;\hat{\psi}\;\in
\;\Gamma(V)
\end{equation}
of the bundle $V$, or rather $V\otimes\C$. The orthogonality
relations of this frame are given by \eqref{eq:normalized_lift}
and \eqref{eq:psi_hat}, namely
\begin{subequations}\label{eq:frame_ortho}
\begin{gather}
\langle \psi,\psi\rangle = \langle \hat{\psi},\hat{\psi}\rangle =0\,,\qquad \langle \psi,\hat{\psi}\rangle =-1\label{eq:frame_a}\\
\langle {\psi},d\psi\rangle =\langle \hat{\psi},d\psi\rangle =\langle d\hat{\psi},\hat{\psi}\rangle =0\,,\label{eq:frame_b}\\
\langle \psi_z,\psi_{z}\rangle =\langle \psi_{\bar{z}},\psi_{\bar{z}}\rangle =0\,,\label{eq:frame_c}\\
\langle \psi_z,\psi_{\bar{z}}\rangle =\tfrac{1}{2}\,.\label{eq:frame_d}
\end{gather}
\end{subequations}
We immediately obtain the fundamental equation of M\"obius
invariant surface geometry, an inhomogeneous Hill's equation,
\begin{equation}\label{eq:c_and_k}
\psi_{z z}+\tfrac{c}{2}\psi=\kappa
\end{equation}
defining the complex valued function $c$ and the section $\kappa$
of $V^{\perp}\otimes\C$. In light of \eqref{eq:Hill_classical} and
\eqref{eq:c} we call $c$ the {\em Schwarzian derivative} of the
immersion $f$ in the coordinate $z$, and denote it by
\begin{equation}\label{eq:Schwarzian_imm}
S_z(f):=c\,.
\end{equation}
We shall see below that the section $\kappa$ can be identified with the
normal bundle valued {\em Hopf differential} of the immersion $f$,
suitably scaled. Note that by construction both, $c$ and $\kappa$,
are M\"obius invariants of the immersion $f$ by given coordinate
$z$. For surfaces in 3--space they form a complete set of
invariants.

If the Hopf differential $\kappa\equiv 0$ is identically zero, the
surface $f$ is totally umbilic and we can view $f$ as a conformal
map taking values in $S^2=\C\P^1$. In this case we recover the
classical Schwarzian discussed in the previous
section~\ref{subsec:classical}.

Just as in the classical case, it is helpful to understand the
behaviour of $c$ and $\kappa$ under coordinate changes. Let $w$ be
another holomorphic coordinate on $M$. Then the new normalized
lift is given by \eqref{eq:lift_change}. Computing
$\tilde{\psi}_{ww}$ and inserting into \eqref{eq:c_and_k} gives
\begin{equation}\label{eq:mess}
2(\tilde{\kappa}-\kappa)z_w^2|w_z|+(\tilde{c}|w_z|-c z_w^2|w_z|+
2z_w^2|w_z|_{zz}+2z_{ww}|w_z|_z)\psi=0\,.
\end{equation}
Taking the $V^{\perp}$ part we see that
\begin{equation}\label{eq:k_transform}
\tilde{\kappa}\tfrac{dw^2}{|dw|}=\kappa\tfrac{dz^2}{|dz|}\,.
\end{equation}
Thus $\kappa\tfrac{dz^2}{|dz|}$ is a globally defined quadratic
differential with values in $LV^\perp\otimes\C$ where $L$ is the real
line bundle $(\bar{K}K)^{-1/2}$ of densities of conformal weight $1$
\cite{Cal98}.

The $\psi$--part of \eqref{eq:mess} gives
\[
\tilde{c}dw^2=(c-2v_{zz} + 2v_z^2)dz^2\,,
\]
where we have put $|w_z|=e^{v}$. From
\eqref{eq:Schwarzian_Euclidean} we therefore deduce the same
transformation rule
\begin{equation}\label{Schwarzian_change}
S_w(f)dw^2=(S_z(f)-S_z(w))dz^2
\end{equation}
for the Schwarzian of a conformal immersion as for the Schwarzian
of a holomorphic function.

We conclude this section by expressing the M\"obius invariants
$S_z(f)$ and $\kappa$ in terms of Euclidean quantities of the
immersion $f:M\to S^n$. If
\[
|df|^2=e^{2u}|dz|^2
\]
is the induced metric of $f$ in the holomorphic coordinate $z$ on
$M$, then
\[
\psi=(1,f)e^{-u}
\]
is the normalized lift \eqref{eq:normalized_lift} of $f$. Let $N_f
M\subset f^\perp$ be the Euclidean normal bundle of $M$.  We have an
isometric isomorphism $N_f S^n\cong V^\perp$:
\begin{equation}
\label{eq:12}
\xi\mapsto \<H,\xi\>(1,f)+(0,\xi)
\end{equation}
and so a section $\hat{\kappa}$ of $N_f S^n\otimes\C$ with
\[
\kappa=\langle H,\hat{\kappa}\rangle (1,f)+(0,\hat{\kappa})\,.
\]
Inserting into \eqref{eq:c_and_k} we obtain
\begin{equation}\label{eq:Hopfdiff}
\hat{\kappa}\,\tfrac{dz^2}{|dz|}=\frac{II^{(2,0)}}{|df|}
\end{equation}
and
\begin{equation}\label{eq:c_Euclidea}
\tfrac{c}{2}\,dz^2=\langle H,II^{(2,0)}\rangle +(u_{zz}-(u_z)^2)dz^2\,.
\end{equation}
Here $II^{(2,0)}$ denote the $(2,0)$-part of the normal bundle
valued second fundamental form of $f$. This shows that $\kappa$,
up to the isomorphism \eqref{eq:12}, is the trace
free part of the second fundamental form, i.e., the normal bundle
valued Hopf differential, scaled by the square root of the induced
metric. In particular, the {\em Willmore energy} of the conformal
immersion $f$ is given by
\begin{equation}\label{eq:Willmore_energy}
W(f)=\int|\kappa|^2\,.
\end{equation}

\subsection{The integrability equations}\label{subsec:integrability}
Unlike the case of holo\-mor\-phic maps, where any Schwarzian can
be integrated, the invariants for a conformal immersion $f:M\to
S^n$ will have to satisfy integrability conditions. Recall that
$f$ gives rise to the M\"obius invariant splitting
\[
M\times \R^{n+1,1}= V\oplus V^{\perp}\,,
\]
where $\P(\mathcal{L}\cap V)$ is the mean curvature sphere
congruence. A choice of holomorphic coordinate $z$ on $M$ yields a
unique normalized lift \eqref{eq:normalized_lift} which, by Hill's
equation \eqref{eq:c_and_k}, defines the Schwarzian derivative
$S_z(f)$ and the Hopf differential $\kappa$.

To compute the integrability conditions on the Schwarzian and the
Hopf differential we work with the orthogonal frame
\begin{equation}\label{eq:frame}
F=(\psi,\psi_z,\psi_{\bar{z}},\hat{\psi},\xi)
\end{equation}
of $V\oplus V^{\perp}$ where $\xi$ denotes a section of
$V^{\perp}$. Clearly, the integrability equations for this frame
are
\[
F_{z\bar{z}}=F_{\bar{z}z}
\]
or, by using reality conditions,
\begin{subequations}\label{eq:int}
\begin{gather}
\psi_{zz\bar{z}}=\psi_{z\bar{z}z}\,,\label{eq:int_a}\\
\text{Im}\,\hat{\psi}_{z\bar{z}}=0\,,\label{eq:int_b}\\
\text{Im}\,\xi_{z\bar{z}}=0\,.\label{eq:int_c}
\end{gather}
\end{subequations}
To evaluate these equations, we will make frequent use of the
formula
\[
\phi=-\langle \phi,\hat{\psi}\rangle \psi-\langle \phi,\psi\rangle \hat{\psi}+
2\langle \phi,\psi_{\bar{z}}\rangle \psi_z+2\langle \phi,\psi_z\rangle \psi_{\bar{z}}\,,
\]
which expresses a section $\phi$ of $V$ in the basis
\eqref{eq:frame_V}. Together with \eqref{eq:c_and_k} we therefore get
\[
\psi_{zz}=-\tfrac{c}{2}\psi+\kappa\,,\qquad\tfrac{c}{2}=\langle \psi_{zz},\hat{\psi}\rangle =-\langle \psi_z,\hat{\psi}_z\rangle \,.
\]
Moreover,
\[
\psi_{z\bar{z}}=q\psi+\tfrac{1}{2}\hat{\psi}\,,\qquad
q=-\langle \psi_{z\bar{z}},\hat{\psi}\rangle =\langle \psi_{z},\hat{\psi}_{\bar{z}}\rangle \,,
\]
for some real valued function $q$. Therefore
\[
\hat{\psi}_z=-c\psi_{\bar{z}}+2q\psi_z+\chi
\]
where $\chi$ is a section of $V^{\perp}\otimes\C\,$. We now
evaluate the first integrability equation \eqref{eq:int_a}
\[
(\kappa-\tfrac{c}{2}\psi)_{\bar{z}}=(q\psi+\tfrac{1}{2}\hat{\psi})_z
\]
which, using the last relation, unravels to
\[
\kappa_{\bar{z}}=2q\psi_z+(\tfrac{1}{2}c_{\bar{z}}+q_z)\psi+\tfrac{1}{2}\chi\,.
\]
To compare $V$ and $V^{\perp}$-parts of this equation, we note
that the $z$-derivative of a section $\xi$ in $V^{\perp}\otimes\C$
decomposes as
\begin{equation}\label{eq:xi_z_decompose}
\xi_{z}=D_z\xi+\langle \xi,\chi\rangle \psi-2\langle \xi,\kappa\rangle \psi_{\bar{z}}\,,
\end{equation}
where $D$ denotes the connection in the bundle $V^{\perp}$.
Therefore,
\[
D_{\bar{z}}\kappa=\tfrac{1}{2}\chi\,,\qquad
q=-\langle \bar{\kappa},\kappa\rangle \,,\qquad
\tfrac{1}{2}c_{\bar{z}}+q_z=\langle \kappa,\bar{\chi}\rangle 
\]
and \eqref{eq:int_a} is equivalent to the {\em conformal Gauss
equation}
\begin{equation*}
\tfrac{1}{2}c_{\bar{z}}=3\langle D_z\bar{\kappa},\kappa\rangle +\langle \bar{\kappa},D_z\kappa\rangle \,.
\end{equation*}
The second integrability condition \eqref{eq:int_b} becomes
\[
0=\Im\hat{\psi}_{z\bar{z}}=\Im(-c_{\bar{z}}\psi_{\bar{z}}-c\psi_{\bar{z}\bar{z}}+2q_{\bar{z}}\psi_{z}+
2q\psi_{z\bar{z}}+\chi_{\bar{z}})\,.
\]
By inserting the appropriate terms and keeping track of reality
conditions, we therefore conclude that \eqref{eq:int_b} is
equivalent to the {\em conformal Codazzi equation}
\begin{equation*}
\Im(D_{\bar{z}}D_{\bar{z}}\kappa+\tfrac{1}{2}\bar{c}\kappa)=0\,.
\end{equation*}
It is easy to see that for the remaining integrability equation
\eqref{eq:int_c} only the $V^{\perp}$-part gives new conditions.
Using \eqref{eq:xi_z_decompose} we immediately obtain the {\em
conformal Ricci equation}
\begin{equation*}
D_{\bar{z}}D_{z}\xi-D_{z}D_{\bar{z}}\xi
-2\langle \xi,\kappa\rangle \bar{\kappa}+2\langle \xi,\bar{\kappa}\rangle \kappa=0
\end{equation*}
for a section $\xi$ of $V^{\perp}$. Note that
\[
R^{D}_{\bar{z}z}\xi=D_{\bar{z}}D_{z}\xi-D_{z}D_{\bar{z}}\xi
\]
is the curvature of the normal bundle $V^{\perp}$.

For further reference we collect the fundamental equations for a
conformal immersion $f:M\to S^n$ in a given holomorphic coordinate
$z$ on $M$: keeping reality conditions in mind the equations of
the frame \eqref{eq:frame} are given by
\begin{equation}\label{eq:moving_frame}
\begin{split}\
\psi_{z\bar{z}}&=-\langle \kappa,\bar{\kappa}\rangle \psi+\tfrac{1}{2}\hat{\psi}\,,\\
\psi_{zz}&=-\tfrac{c}{2}\psi+\kappa\,,\\
\hat{\psi}_z&=-2\langle \kappa,\bar{\kappa}\rangle \psi_z -c\psi_{\bar{z}}
+2D_{\bar{z}}\kappa\,,\\
\xi_z&=D_z\xi+2\langle \xi,D_{\bar{z}}\kappa\rangle \psi-2\langle \xi,\kappa\rangle \psi_{\bar{z}}\,.
\end{split}
\end{equation}
The resulting integrability conditions, the conformal Gauss,
Codazzi and Ricci equations, are
\begin{subequations}\label{eq:fundamental}
\begin{gather}
 \tfrac{1}{2}c_{\bar{z}}=3\langle D_z\bar{\kappa},\kappa\rangle +\langle \bar{\kappa},D_z\kappa\rangle \,,\label{eq:Gauss}\\
 \Im(D_{\bar{z}}D_{\bar{z}}\kappa+\tfrac{1}{2}\bar{c}\kappa)=0\,,\label{eq:Codazzi}\\
 R^{D}_{\bar{z}z}\xi=D_{\bar{z}}D_{z}\xi-D_{z}D_{\bar{z}}\xi=
2\langle \xi,\kappa\rangle \bar{\kappa}-2\langle
 \xi,\bar{\kappa}\rangle \kappa\,.\label{eq:Ricci}
\end{gather}
\end{subequations}
These equations can be given an invariant co-ordinate free meaning
using the technology of Cartan connections \cite{BurCal}.
However this is not our aim here, where we intend mostly to focus on classical
surface theory in 3-space. In that case the conformal Ricci equation
is vacuous and we get the following fundamental theorem of conformal
surface theory:
\begin{theorem}\label{thm:fundamental}
Let $M$ be a Riemann surface and let $L=(K\bar{K})^{-1/2}$ denote
the $1$-density bundle, where $K$ is the canonical bundle of $M$.

If $f_k:M\to S^3$ are two conformal immersions inducing the same
Hopf differentials and Schwarzians, then there is a M\"obius
transformation $T:S^3\to S^3$ with $Tf_1=f_2$.

Conversely, let $\kappa\tfrac{dz^2}{|dz|}\in \Gamma(LK^2)$ be a quadratic
differential with values in $L$ and $c$ be
a Schwarzian derivative, i.e., $c$ transforms according to
\eqref{eq:Schwarzian_Euclidean}. If $\kappa$ and $c$ satisfy the
conformal Gauss and Codazzi equations then there exists a
conformal immersion $f:M\to S^3$ with M\"obius holonomy, Hopf
differential $\kappa$, and Schwarzian derivative $c$.
\end{theorem}
\subsection{Applications to old and new results in conformal surface theory}
At this stage it is worthwhile to pause and apply our conformally
invariant setup to discuss some results of conformal surface
theory. We content ourselves by showing the basic ideas and will
work locally, leaving the global discussion as an exercise to the
interested Reader.

The simplest question one can ask in any geometry is to
characterize the homogeneous surfaces, i.e., surfaces which are
2-parameter orbits of the symmetry group under consideration. In
M\"obius geometry such surfaces are called the {\em cyclides of
Dupin}. Clearly, such a surface must have Hopf differential and
Schwarzian constant in the coordinates given by the 2-parameter
group. Due to the transformation rule \eqref{eq:Schwarzian_imm}
the constancy of the Schwarzian has the invariant meaning that the
Schwarzian is holomorphic (this amounts to the flatness of a certain
normal Cartan connection on $M$ \cite{BurCal}).
\begin{theorem}\label{thm:Dupin}
Let $f:M\to S^3$ be a conformal immersion, not contained in a
2-sphere, whose Schwarzian is holomorphic. Then $f$ is a cyclide
of Dupin.
\end{theorem}
\begin{proof}
It suffices to show that $\kappa$ and $S_w(f)$ are constant in an
appropriate holomorphic coordinate $w$. Since $c$ is holomorphic
the conformal Gauss equation \eqref{eq:Gauss} reads
\[
(\kappa\bar{\kappa}^3)_z=0\,,
\]
expressing the fact that $\bar{\kappa}\kappa^3$ is holomorphic. By
our assumption $f$ is not contained in a 2-sphere so that $\kappa$
is non-zero. Therefore
\[
dw=(\bar{\kappa}\kappa^3)^{1/4}dz
\]
defines a new holomorphic coordinate $w$ on $M$. Due to the
transformation rule \eqref{eq:k_transform} for sections of
$LK^2$ we see that
\[
\kappa=(w_z)^{3/2} \bar{w_z}^{-1/2}\tilde{\kappa}\,.
\]
Combining these two equations gives $\tilde{\kappa}=1$ and the
conformal Codazzi equation \eqref{eq:Codazzi} forces $\tilde{c}$
to be imaginary, and thus constant.
\end{proof}
A further basic question to ask is how much of the invariants of
an immersed surface are actually needed to determine the surface
up to symmetry. In the Euclidean setup this is the Bonnet problem:
generically the induced metric and the mean curvature determine a
surface up to rigid motions. In M\"obius geometry one expects that
generically the Hopf differential should determine a surface up to
M\"obius transformations.

However, let us recall the notion of {\em isothermic} surfaces
\cite{Bia05,Cal03,Dar99}. Such
surfaces are conformally parameterized by their curvature lines
away from umbilic points. Put differently, away from umbilics
there are holomorphic coordinates in which the Hopf differential
$\kappa$ is {\em real valued}. In this case the conformal Gauss
and Codazzi equations simplify to
\begin{gather*}
c_{\bar{z}}=4(\kappa^2)_z\,,\\
\Im\,(\kappa_{\bar{z}\bar{z}}+\tfrac{1}{2}\bar{c}\kappa)=0\,.
\end{gather*}
Assuming $\kappa$ non-zero the two equations combine to Calapso's
equation \cite{Cal03} 
\[
\triangle\tfrac{\kappa_{xy}}{\kappa}+8(\kappa^2)_{xy}=0
\]
where $z=x+iy$. The conformal Gauss and Codazzi equations are
invariant under deformations of the Schwarzian by
\[
c_{r}=c+r
\]
where $r\in\R$ is a real parameter, the spectral parameter of
algebro-geometric integrable system theory.  By
theorem~\ref{thm:fundamental} we thus see that an isothermic surface
comes in an associated family parameterized by $\R$: these are the
$T$-transforms of Calapso and Bianchi \cite{Bia05A,Cal03}. Note that
since all $c_{r}$ are distinct the surfaces in the associated
family are non-congruent with respect to the M\"obius group.
Therefore the Hopf differential alone does not determine isothermic
surfaces. But these surfaces are the only exceptions:
\begin{theorem}\label{thm:fundamental_noniso}
Let $f_k:M\to S^3$ be two non-congruent conformal immersions inducing
the same Hopf differentials. Then $f_1$ and $f_2$ are isothermic
surfaces in the same associated family.
\end{theorem}
\begin{proof}
Since $\kappa=\kappa_k$ the conformal Codazzi equation implies
that the difference of the Schwarzians
\[
c dz^2:=c_1dz^2-c_2dz^2\in H^0(K^2)
\]
is a holomorphic quadratic differential. Since the $f_k$ are
non-congruent, $c$ is not identically zero, and so we can choose
holomorphic coordinates such that
\[
c=1.
\]
The conformal Codazzi equation now implies
\[
\Im\,(\kappa)=0\,,
\]
so that both surfaces are isothermic and in the same associated family.
\end{proof}
We conclude this section with a brief discussion of (constrained)
Willmore surfaces: it can be checked that a conformal immersion
$f:M\to S^n$ is Willmore if and only if
\[
\Re\,(D_{\bar{z}}D_{\bar{z}}\kappa+\tfrac{1}{2}\bar{c}\kappa)=0\,.
\]
Therefore the integrability conditions for a Willmore surface are
\begin{gather*}
D_{\bar{z}}D_{\bar{z}}\kappa+\tfrac{1}{2}\bar{c}\kappa=0\,,\\
\tfrac{1}{2}c_{\bar{z}}=3\langle D_z \bar{\kappa},\kappa\rangle +\langle \bar{\kappa},D_z\kappa\rangle \,,\\
R^{D}_{\bar{z}z}\xi=
2\langle \xi,\kappa\rangle \bar{\kappa}-2\langle \xi,\bar{\kappa}\rangle \kappa \,.
\end{gather*}
More generally, as we shall see below, a surface is \emph{constrained
Willmore}---that is, it extremizes the Willmore functional with
respect to variations through conformal immersions---if
\begin{equation}
\label{eq:15}
D_{\bar{z}}D_{\bar{z}}\kappa+\tfrac{1}{2}\bar{c}\kappa=\Re(\bar{q}\kappa)
\end{equation}
for some holomorphic quadratic differential $q dz^2\in H^0(K^2)$.
These systems have the symmetry
\begin{equation}
\label{eq:13}
\kappa_{\lambda}=\lambda\kappa,\qquad 
c_\lambda=c+(\lambda^2-1)q\qquad
q_\lambda=\lambda^2 q 
\end{equation}
for unitary $\lambda\in S^1$, which describes the associated family
of (constrained) Willmore surfaces in our setup. This spectral
deformation is particularly simple for Willmore surfaces ($q$=0):
here the Schwarzian is fixed and the scaled Hopf differential
$\kappa\tfrac{dz^2}{|dz|}\in\Gamma(LK^2)$ is rotated by a phase.  All
this is reminiscent of the theory of constant mean curvature surfaces
in $3$-space where the associated family is obtained by rotating the
Euclidean Hopf differential by a phase and fixing the induced
metric.  As we shall see, this similarity is not a coincidence.

A classical result of Thomsen \cite{Tho23} characterizes isothermic
Willmore surfaces in $3$-space as minimal surfaces in some space-form
(conic section).  More generally, Richter \cite{Ric97} proves that
constant mean curvature surfaces in space forms are isothermic and
constrained Willmore.  These results are easy to see in our setup:
from \eqref{eq:1} we see that a surface has mean curvature $H$ in the
conic section $S_{v_0}$ exactly when $H=-\<v_0,N\>$ for $N$ a unit
length section of $V^\perp$ (in particular, a surface is minimal if
$V$ contains a constant section $v_0$).  Moreover, the metric
induced from $S_{v_0}$ is given by
\begin{equation}
\label{eq:2}
e^{2u}|dz|^2=\frac{|dz|^2}{\<\psi,v_0\>^2}.
\end{equation}
Suppose then that $f$ has constant mean curvature $H$ in $S_{v_0}$
with induced metric $e^{2u} |dz|^2$.  We write
\begin{equation}
\label{eq:3}
v_0=a\psi+b\psi_z+\bar{b}\psi_{\bar{z}}+\hat{a}\hat{\psi}-HN
\end{equation}
for real valued functions $a$, $\hat{a}$ and a complex valued
function $b$.  Writing out $v_{0,z}$ in terms of our frame and
using \eqref{eq:moving_frame}, we see that the constancy of $v_0$
amounts to
\begin{subequations}
\label{eq:4}
\begin{gather}
\label{eq:5}
b\kappa+2\hat{a}\kappa_{\bar{z}}=0\\
\label{eq:6}
\hat{a}_z+\tfrac{1}{2}\bar{b}=0\\
\label{eq:7}
\bar b_z-c\hat a +2H\kappa=0\\
\label{eq:8}
a+b_z-2\hat a |\kappa|^2=0\\
\label{eq:9}
a_z-\tfrac{c}{2}b-\bar b|\kappa|^2-2H\kappa_{\bar z}=0
\end{gather}
\end{subequations}
{}From \eqref{eq:2} we have $\hat a=e^{-u}$ and then \eqref{eq:5} and
\eqref{eq:6} yield
\[
e^{-2u}(e^u\kappa)_{\bar z}=0
\]
so that the Hopf differential $e^u\kappa dz^2$ is holomorphic.  Thus,
away from umbilics, we can choose the holomorphic coordinate so that
$e^u\kappa=1$ or, equivalently, that $\hat a=\kappa$.  In particular,
$\kappa$ is real and the surface is isothermic.  Now \eqref{eq:6}
gives $b=-2\kappa_{\bar z}$ while \eqref{eq:8} gives
$a=2(\kappa^3+k_{z\bar z})$.  This leaves \eqref{eq:7} which reads
\[
\kappa_{\bar z\bar z} + \tfrac{\bar c}{2}\kappa=H\kappa
\]
so that $\psi$ is Willmore if $H=0$ and constrained Willmore
otherwise, while, in view of the Gauss equation, \eqref{eq:9} is an
identity.

Conversely, given a surface with $\kappa$ real and
\begin{equation}
\label{eq:10}
\kappa_{\bar z\bar z} + \tfrac{\bar c}{2}\kappa=H\kappa
\end{equation}
for some constant $H$ we define $v_0$ by \eqref{eq:3} with the same
choice of coefficients as above and conclude from \eqref{eq:4} that
$v_0$ is constant whence $\psi$ gives a surface of constant mean
curvature $H$ in the space form $S_{v_0}$ of constant curvature $K$
given by
\begin{equation}
\label{eq:11}
K=-\<v_0,v_0\>=-H^2+4\kappa^4+4\kappa\kappa_{z\bar z}-4|\kappa_z|^2
\end{equation}
(see \cite{MusNic01} for a recent analysis of the soliton theory of this
equation).

The isothermic spectral deformation $c_r=c+r$, $r\in\R$, preserves
the class of constant mean curvature surfaces: in view of
\eqref{eq:10}, we have
\[
\kappa_{\bar z\bar z} + \tfrac{\bar{c_r}}{2}\kappa=
(H+\tfrac{r}{2})\kappa
\]
so that the corresponding surface $\psi_r$ has constant mean
curvature $H_r=H+\tfrac{r}{2}$ is a space form of
curvature $K_r$ where, in view of \eqref{eq:11},
$K_r+H_r^2$ is independent of $r$.  Moreover, all
these surfaces are isometric with metric $\kappa^{-2}|dz|^2$.  Thus our
spectral deformation recovers the family of constant mean curvature
surfaces related by the Lawson correspondence \cite{Law70}.

Moreover, the spectral symmetry \eqref{eq:13} of
constrained Willmore surfaces preserves the isothermic condition: if
$\kappa$ is real for a coordinate $z$ then, by
\eqref{eq:k_transform}, $\kappa_\lambda$ is real in the rotated
coordinate $w=\sqrt{\lambda}z$.

For constant mean curvature surfaces, these two deformations can be
combined to give an action of $\C\setminus\{0\}$ on such surfaces.  Indeed,
for a surface of constant mean curvature $H$ in a space form of
curvature $K$ , we have seen that there is a coordinate $z$ for which
$\kappa$ is real and $q=H$.  For $\lambda\in\C\setminus\{0\}$, set
\begin{equation*}
\kappa_\lambda=\tfrac\lambda{|\lambda|}\kappa\qquad
c_\lambda=c+2(\lambda-1)q\qquad
q_\lambda=\lambda q.
\end{equation*}
The fact that $\bar{q}\kappa$ is real ensures that each
$\kappa_\lambda,c_\lambda, q_\lambda$ solves \eqref{eq:fundamental}
and \eqref{eq:15}.  Further, in the rotated coordinate
$w=\sqrt{\frac\lambda{|\lambda|}}z$, $\kappa_\lambda$ is real (in
fact $\kappa_\lambda dz^2=\kappa dw^2$) while $q_\lambda dz^2 =
|\lambda|H dw^2$ so that, with respect to $w$, we have a solution of
\eqref{eq:10} with $H_\lambda=|\lambda| H$.  From \eqref{eq:11}, we
conclude that the surface $\psi_\lambda$ has constant mean curvature
$|\lambda|H$ and metric $\kappa^{-2}|dz|^2$ in a space form of
curvature $K_\lambda=K+(1-|\lambda|^2)H$.  In particular, for
$\lambda\in S^1$, we have a family of isometric surfaces in a fixed
space form with fixed $H$ and rotated Hopf differential---this is the
classical associated family.  For minimal surfaces, the action of
$\R\setminus\{0\}$ is trivial but when $H\neq 0$ we get (part of) the Lawson
deformation so that these two well-known symmetries are unified in a
single $\C\setminus\{0\}$ action. 

For $M$ a torus, any holomorphic quadratic differential is of the
form $qdz^2$ with $q$ constant so that we conclude with Richter \cite{Ric97}
that \emph{all} constrained Willmore isothermic tori have constant
mean curvature in some space form.  However, without this global
hypothesis, there remains the possibility of isothermic constrained
Willmore surfaces that are not spectral deformations of isothermic
Willmore surfaces.

One might ask if a version of Thomsen's result holds in $n$-space.
We restrict attention to isothermic Willmore surfaces of which
obvious examples, in addition to the above mentioned minimal
surfaces, are cylinders in $\R^4$ over elastic curves in $\R^3$,
cones over elastic curves in $S^3$ and rotational surfaces about the
plane at infinity of elastic curves in hyperbolic $3$-space $H^3$.
These surfaces are clearly isothermic and, by the principle of
symmetric criticality, Willmore. 
\begin{theorem}
Let $f:M\to S^n$ be an isothermic Willmore surface  with $n\geq
3$. Then either $f$ is minimal in $3$-space, or an
isothermic Willmore surfaces in $4$-space described by an ODE:
the equation describing elastica or a Painleve-type equation
for the Hopf differential with radial symmetry.
\end{theorem}
We begin by collecting the integrability equations for isothermic
Willmore surfaces: since $\kappa=\bar{\kappa}$ is real, Ricci's
equation immediately implies that the normal bundle is flat. We
can therefore assume $V^{\perp}=M\times \R^{n-2}$ and $D$ is just
directional derivative. The two remaining integrability
conditions, the Gauss and Codazzi equations, become
\begin{gather*}
\kappa_{zz}+\tfrac{c}{2}\kappa=0\,,\\
c_{\bar{z}}=4\langle \kappa,\kappa\rangle _z\,,
\end{gather*}
for the $\C$-valued Schwarzian $c$ and the $\R^{n-2}$-valued Hopf
differential $\kappa$. Differentiating the Codazzi equation twice
with respect to $\bar{z}$ and keeping reality conditions and the
Gauss equation in mind, we obtain
\begin{equation}\label{eq:gradient}
\text{Im}\,\langle \kappa,\kappa\rangle _z\kappa_{\bar{z}}=0\,.
\end{equation}
First note that this is automatically satisfied if $n=3$, in which
case we already know that $f$ is minimal. From now on we assume
that $n\geq 4$. If $\kappa$ is a nonzero constant, then Gauss'
equation implies that $c=0$ and by using the frame equations
\eqref{eq:moving_frame}, one can check that $f$ is the Clifford
torus. The case $\kappa=0$ results, as we already have discussed,
in a holomorphic map $f$ into the $2$-sphere. Therefore we may
assume that $\kappa$ is non-constant. If $\kappa$ has constant length
$\langle\kappa,\kappa\rangle\neq 0$
then $\langle\kappa_z,\kappa\rangle=0$, and differentiating
yields 
\[ 
0= \langle\kappa_{zz},\kappa\rangle+ \langle\kappa_z,\kappa_z\rangle = - \tfrac{c}{2} 
\langle\kappa,\kappa\rangle
+\langle\kappa_z,\kappa_z\rangle\,. 
\] 
{}From the Gauss and Codazzi equations we see that $c$ is holomorphic and that
$\langle\kappa_z,\kappa_z\rangle$ is
anti-holomorphic. Therefor both are constant.
After a rotation of the coordinate we may assume that $c$ is a real
constant. Then
$\kappa$ is contained in a rotated $\R^{n-2}$ in $\C^{n-2}$ and
has constant length.
Decomposing the complex Codazzi equation according to 
the splitting $\C^{n-2}=\R^{n-2}\oplus i\R^{n-2}$, we get  
\[
\kappa_{xx} - \kappa_{yy} + 2 c \kappa = 0\,,\qquad\text{and}\qquad 
\kappa_{xy} =0\,.
\]
The latter implies $\kappa(x,y) =
k_1(x)+k_2(y)$ which, inserted into the former, gives
\[ 
f''(x) + 2 c f(x) = g''(y) - 2c g(y)=\text{constant}\,. 
\] 
So $k_1$ and $k_2$ solve linear ordinary differential equations
which can be explicitly integrated. If $c=0$ only constant $\kappa$'s are solutions 
of constant length. If $c>0$ then any constant length solution is of the form 
\[ 
\kappa(x,y) = \cos(\sqrt{2c}\, x) \kappa_1 + \sin(\sqrt{2c}\,x) \kappa_2 
\]
for $\kappa_1$, $\kappa_2\in \R^{n-2}$ two orthogonal vectors of
same length. The corresponding surfaces are the 
rotational surfaces, cylinders and cones over free elastic helixes in the
respective $4$-dimensional space forms. 

We now may assume that the length of $\kappa$ is non-constant. In this case 
\eqref{eq:gradient} is equivalent to
\[
\kappa=k(g)\qquad\text{with}\qquad g=\langle \kappa,\kappa\rangle \,,
\]
for some $\R^{n-2}$-valued function $k$ of one variable. Inserting
into the Codazzi equation yields the second order differential equation
\[
g_z^2 k''+g_{zz} k'+\tfrac{c}{2}k=0
\]
for the function $k$. This shows that the image of $\kappa$ is
contained in an at most $2$-dimensional parallel subbundle of
$E\subset V^{\perp}$. Using the fourth frame equation in
\eqref{eq:moving_frame}, we deduce that the complementary
subbundle $E^{\perp}\subset V^{\perp}$ is a constant subspace in
the ambient space. We therefore conclude that $f$ takes values in
$S^4$.

Assuming from now on that our isothermic Willmore surface $f$ is
fully contained in $4$-space, its trivial normal bundle has a
well-defined $90$-degree rotation $J$. Since $\langle \kappa,J\kappa\rangle =0$,
the Codazzi equation gives
\[
\langle \kappa_z,J\kappa\rangle _z=0\,.
\]
{}From this we deduce that
\[
\langle \kappa_z,J\kappa\rangle ={\bar{h}}^{-1}\,,\qquad h_{\bar{z}}=0
\]
for a holomorphic function $h$. Note that $h$ is nonzero,
otherwise $\kappa$ would have values in a $1$-dimensional parallel
subbundle of the normal space and thus $f$ would not be full in
$4$-space. Since $\kappa=k(g)$ we get
\[
\tfrac{h}{|h|^2}=\langle \kappa_z,J\kappa\rangle =\langle k'(g),Jk(g)\rangle g_z
\]
implying that $h$ is a real multiple of $g_z=\langle \kappa,\kappa\rangle _z$.
Multiplying \eqref{eq:gradient} by this real multiple gives
\[
\text{Im}\,h\kappa_{\bar{z}}=0\,.
\]
But this is equivalent to $\kappa$ depending only on
$\text{Re}\,w$, where
\[
w_z=h
\]
is a new holomorphic coordinate. We may thus assume that
\[
\kappa=k(w+\bar{w})
\]
for some $\R^2$-valued function of one variable, which we again
call $k$. Inserting back into the Codazzi equation yields
\[
h^2 k''+h_z k'+\tfrac{c}{2}k=0\qquad\text{or}\qquad
k''+\tfrac{h_z}{h^2}k'+\tfrac{c}{2h^2}k=0\,.
\]
Now $f$ is assumed to be full in $4$-space so that $k$ and $k'$
are linearly independent. Since $k$ is real, we conclude that
\[
h_z=\alpha h^2\,,\qquad\alpha\in\R
\]
and then that
\begin{equation}\label{eq:Schwarz}
\frac{c}{2h^2}=\rho(w+\bar{w})\,.
\end{equation}
The first relation tells us that $h=\tfrac{1}{az+b}$ is a special
fractional linear coordinate change. Taking $\bar{z}$ derivative
of the  second relation and applying the Gauss equation unravels to
\[
\rho'=\frac{2\langle k,k\rangle '}{|h|^2}\,.
\]
We have two cases to discuss: first, if $h$ is constant
the last relation can be integrated to
\[
\rho=2\frac{\langle k,k\rangle +\gamma}{|h|^2}\,,\qquad
\gamma\in\R\,,\,\,\,h\in\C\,.
\]
This yields the formula
\begin{equation}
c=2\frac{h}{\bar{h}}(\langle k,k\rangle +\gamma)
\end{equation}
for the Schwarzian of our surface $f$. Inserting back into the Codazzi
equation gives the second order ODE for spatial elastica
\begin{equation}\label{eq:k}
k''+\frac{h}{\bar{h}}(\langle k,k\rangle +\gamma)k=0
\end{equation}
the Hopf differential $\kappa=k(hz+\bar{h}\bar{z})$ has to satisfy.

In the second case the fractional linear coordinate change is,
after a translation and real scale, of the form
\[
h(z)=\frac{1}{z}\,.
\]
Then $w+\bar{w}=2\ln |z|$ and 
\[
\rho'=2\langle k,k\rangle |z|^2\,,\qquad k''-k+\rho k=0\,,
\]
is and ODE for the Hopf differential $k$ depending on the radius
only. This case is similar to the intrinsic rotational constant
mean curvature surfaces \cite{BobIts95,Smy93,TimPinFer94}.
 
Finally, solving the respective ODE's
determines the Hopf differential and the Schwarzian is computed from
\eqref{eq:Schwarz}. Inserting the invariants into the frame equations
\eqref{eq:moving_frame} allows to determine the surface $f$.
We leave the detailed calculations to the interested Reader.

\section{The Novikov--Veselov and Davey--Stewartson flows}
Over the last 15 years a number of papers have been written on
{\em the soliton theory of surfaces}. In all cases only special
surface classes, like constant (mean) curvature, isothermic or
Willmore surfaces, were considered. In case the underlying surface
is a 2-torus it has been shown that those surfaces allow an
infinite hierarchy of commuting flows. Surfaces stationary under a
finite number of flows could be explicitly constructed by noticing
that those flows are given by the linear flows on the Jacobian of
some auxiliary Riemann surface, the so-called {\em spectral
curve}.

More recently, Konopelchenko and his collaborators
\cite{Kon00,KonLan00,KonTai95A,Tai97} have shown how to construct
hierarchies of commuting flows on the space of \emph{all} conformally
immersed tori in $3$ and $4$-dimensional spheres based on the
modified Novikov--Veselov and Davey--Stewartson hierarchies.
 
Following our construction of the KdV-hierarchy, we are now going to
give a manifestly conformally invariant construction of such flows
and show how the special surface classes mentioned above are
invariant under them.

\subsection{Deformations of conformal immersions}
Let $M=T^2$ be a 2-torus and $X\in H^0(TM)$ a holomorphic vector
field so that $\tfrac{\del}{\del z}=\tfrac{1}{2}(X-iJX)$ for a
holomorphic coordinate $z$ on the universal cover $\R^2$. Given a
deformation $f(t):M\to S^n$ of conformal immersions with $f=f(0)$,
we get the unique normalized lifts
\[
\psi(t):M\to \mathcal{L}_{+}
\]
from \eqref{eq:normalized_lift} satisfying
\[
|d\psi(t)|^2=|dz|^2\,.
\]
Therefore we have the canonical frame \eqref{eq:frame_V} of $V(t)$
for all times $t$. The deformation of this frame is determined by
the variation
\begin{equation}\label{eq:Psi_dot}
\psi_t=a\psi+b\psi_z+\bar{b}\psi_{\bar z}+\sigma
\end{equation}
of the normalized lift $\psi$. Here $\sigma(t)\in
\Gamma(V(t)^{\perp})$ is the normal variation, the complex valued
function $b(t):M\to \C$ determines the tangential and
$a(t):M\to\R$ the radial variation. Since $\psi$ is null, there is
no $\hat{\psi}$ component. For the deformation $f(t)$ to remain
conformal, i.e.
\begin{equation}\label{eq:conf_def}
|d\psi|^2_t=0\,,
\end{equation}
the radial and tangential deformations $a$ and $b$ have to be
determined by the normal deformation $\sigma$. To compute those
relations, we need
\begin{equation}\label{eq:psi_z_t}
\begin{split}
(\psi_z)_t&=(\psi_t)_z\\
&=(a_z-\tfrac{c}{2}b-\bar{b}|\kappa|^2+2\langle \sigma,D_{\bar
z}\kappa\rangle )\psi+(a+b_z)\psi_z+(\bar{b}_z-2\langle \sigma,\kappa\rangle )\psi_{\bar
z} +\tfrac{\bar{b}}{2}\hat{\psi}+(b\kappa+D_z\sigma)
\end{split}
\end{equation}
where we have used the frame equations \eqref{eq:moving_frame}.
Inserting into \eqref{eq:conf_def} gives
\begin{equation}\label{eq:conf_def_eq}
\bar{b}_z=2\langle \sigma,\kappa\rangle \qquad\text{and} \qquad a=-\Re\,b_z\,,
\end{equation}
of which the first is a $\delbar$-problem solvable over a torus if
and only if the right hand side has vanishing integral $\int_M
\langle \sigma,\kappa\rangle =0$.

Using \eqref{eq:conf_def_eq} we can eliminate $a$ in
\eqref{eq:psi_z_t} and rewrite
\begin{equation}\label{eq:psi_z_dot}
(\psi_z)_t=(-\tfrac{1}{2}b_{zz}-\tfrac{c}{2}b-\bar{b}|\kappa|^2+\langle \sigma,D_{\bar
z}\kappa\rangle -\langle D_{\bar z}\sigma,\kappa\rangle )\psi+(\Im\,b_z)
\psi_z+\tfrac{\bar{b}}{2}\hat{\psi}+(b\kappa+D_z\sigma)\,.
\end{equation}
(Here $\Im b_z=\tfrac{1}{2}(b_z-\bar{b_z})$ and so is pure
imaginary.)

We now can calculate the remaining deformations. First note that
\begin{equation}\label{eq:psi_hat_dot}
\hat{\psi}_t=(\Re\,b_z)\hat{\psi}-
4\Re\,(\langle \hat{\psi},\psi_{\bar{z}t}\rangle \psi_z)+\tau
\end{equation}
where, using \eqref{eq:psi_z_dot}, we have
\[
2\langle \hat{\psi},\psi_{zt}\rangle =b_{zz}+cb+2\bar{b}|\kappa|^2-2\langle \sigma,D_{\bar
z}\kappa\rangle +2\langle D_{\bar z}\sigma,\kappa\rangle \,.
\]
Moreover, the normal variation $\tau$ is obtained from
\eqref{eq:moving_frame},
\[
\tau=(\hat{\psi}_t)^{\perp}=2((\psi_{zt\bar{z}})^{\perp}+|\kappa|^2\sigma)\,,
\]
which unravels to
\[
\tfrac{\tau}{2}=D_{\bar z}D_z\sigma
+2\langle \sigma,\bar{\kappa}\rangle \kappa+|\kappa|^2\sigma+2\Re\,(bD_{\bar
z}\kappa)\,.
\]
Note that the latter is indeed real due to the Ricci equation
\eqref{eq:Ricci}. Finally, if $\xi$ is normal then its tangential
deformation calculates to
\begin{equation}\label{eq:xi_dot_tan}
(\xi_t)^{T}=\langle \xi,\tau\rangle \psi+\langle \xi,\sigma\rangle \hat{\psi}-
4\Re\,(\langle \xi,b\kappa+D_z\sigma\rangle \psi_{\bar z})\,.
\end{equation}
To calculate the deformations of the Hopf differential $\kappa$
and the Schwarzian $c$, we use Hill's equation \eqref{eq:c_and_k}
and obtain
\[
\kappa_t=\psi_{tzz}+\tfrac{c_t}{2}\psi +\tfrac{c}{2}\psi_t\,.
\]
Taking the $V^{\perp}$ part, together with \eqref{eq:psi_z_dot} and
\eqref{eq:moving_frame}, gives
\begin{equation}\label{eq:kappa_dot}
D_t\kappa=D_z D_z\sigma
+\tfrac{c}{2}\sigma+(\text{Im}b_z+b_z)\kappa+bD_z\kappa+\bar{b}D_{\bar{z}}\kappa\,.
\end{equation}
Similarly, using \eqref{eq:psi_hat_dot} and taking an inner product
with $\hat\psi$ we get
\begin{equation}\label{eq:C_dot}
\begin{aligned}
\tfrac{1}{2}c_t=&\tfrac{1}{2}b_{zzz}+cb_z+\tfrac{1}{2}(bc_z+\bar{b}c_{\bar{z}})+
8\langle \sigma,\bar{\kappa} \rangle\langle \kappa,\kappa\rangle +\\&
3\langle D_{\bar{z}}D_z\sigma,\kappa\rangle -\langle \sigma,D_{\bar{z}}D_z\kappa\rangle -
3\langle D_z\sigma,D_{\bar{z}}\kappa\rangle +\langle D_{\bar{z}}\sigma,D_z\kappa\rangle \,.
\end{aligned}
\end{equation}
For later use, we need to record the deformation of the Willmore
energy \eqref{eq:Willmore_energy}. Up to exact forms
\[
\tfrac{1}{2}|\kappa|^2_t=\text{Re}\langle D_t\kappa,\bar{\kappa}\rangle \equiv
\text{Re}\langle D_z D_z\sigma+\tfrac{c}{2}\sigma,\bar{\kappa}\rangle \equiv
\langle \sigma,D_z D_z\bar{\kappa}+\tfrac{c}{2}\bar{\kappa}\rangle \,,
\]
where we used \eqref{eq:kappa_dot} together with the product rule.
Note that due to Codazzi's equation \eqref{eq:Codazzi} the last
expression is real. Therefore the Willmore energy deforms by
\begin{equation}\label{eq:W_dot}
W_t=2\int\langle \sigma,D_z D_z\bar{\kappa}+\tfrac{c}{2}\bar{\kappa}\rangle \,.
\end{equation}
Note that $W_t=0$ for compactly supported variations if and only
if
\[
D_{\bar{z}}D_{\bar{z}}\kappa+\tfrac{1}{2}\bar{c}\kappa=\text{Re}(\bar{q}\kappa)
\]
for some holomorphic quadratic differential $q dz^2$. This
follows from the constraint $\int\langle \sigma,\kappa\rangle =0$ on the normal
variations \eqref{eq:conf_def_eq}, expressing the fact that we
deform via {\em conformal} immersions $f(t):M\to S^n$ rather then
{\em all} immersions. The resulting critical points are the {\em
constrained Willmore} immersions. Every Willmore immersion,
$q=0$, clearly is constrained Willmore.
\subsection{The KdV flows revisited}
As a warm up, we make contact with our discussion of the KdV flows
earlier on. We already have indicated that holomorphic maps into
$\C\P^1$ are characterized by vanishing Hopf differential $\kappa
\equiv 0$. In this case the $\delbar$-problem
\eqref{eq:conf_def_eq} has the holomorphic solution $b=c$.
Inserting into the variation of the Schwarzian \eqref{eq:C_dot},
we immediately obtain the KdV equation
\[
c_t=c_{zzz}+3c_z c\,.
\]
\subsection{The Novikov--Veselov flows}\label{subsec:NV_flow}
Adapting our recursive scheme in the construction of the KdV
hierarchy in Section~\ref{subsec:KdV_hierarchy} to deformations of
conformal immersion, we outline the first steps in the construction
of a commuting hierarchy of odd-order flows on conformally immersed
tori. For these flows to restrict to the known flows on say Willmore
surfaces, they should at least preserve the Willmore energy
\eqref{eq:Willmore_energy}. Contrary to the existing literature
\cite{Tai97}, where a hierarchy of flows for surfaces in $3$-space is
described in terms of Euclidean data, our flows are M\"obius
invariant from the onset.

Any conformal deformation $f(t):M\to S^n$ is determined, up to the
solution of a $\delbar$-problem \eqref{eq:conf_def_eq}, by the
normal variation \eqref{eq:Psi_dot}. The most naive flow, which we
will label as the first flow, has no normal variation, i.e.,
$\sigma=0$, and thus $b=1$. Then \eqref{eq:kappa_dot},
\eqref{eq:C_dot} give
\[
D_{t_1}\kappa=(D_z+D_{\bar{z}})\kappa\qquad\text{and}\qquad
c_{t_1}=\tfrac{1}{2}(c_z+c_{\bar{z}})\,,
\]
which is just translational flow in $X$ direction.

Our previous scheme---the new normal deformation is given by the
old variation of the invariants---suggests $\sigma=D_{t_1}\kappa$
for the construction of the next flow up. This is certainly
reasonable, since both, $\kappa$ and $\sigma$, are normal bundle
valued. Since $\sigma$ has to be {\em real} and the resulting flow
should preserve the Willmore energy \eqref{eq:Willmore_energy},
\begin{equation}\label{eq:sigma_NV}
\sigma=\text{Re}\,D_z\kappa
\end{equation}
turns out to be the right choice for the third flow. Let us first
check that we can solve the $\delbar$-problem
$\bar{b}_z=2\langle \sigma,\kappa\rangle $: modulo exact forms, we obtain
\[
2\langle \text{Re}\,D_z\kappa,\kappa\rangle \equiv
\langle D_{\bar{z}}\bar{\kappa},\kappa\rangle \equiv
-\langle \bar{\kappa},D_{\bar{z}}\kappa\rangle \equiv 0\,,
\]
where the last identity uses the Gauss equation \eqref{eq:Gauss}.
Thus we can solve for $b$, and hence for $a$, in
\eqref{eq:conf_def_eq} to get the variation \eqref{eq:Psi_dot} of
$\psi$ respectively $f$. The solution $b$ of the $\delbar$-problem
on the torus $M$ is determined only up to a constant for each time
$t$. To obtain a well defined flow in $t$, we normalize $b$ so
that $\int b=0$.

To see whether the Willmore energy is preserved, we calculate
\eqref{eq:W_dot} the expression
\begin{equation}\label{eq:NV_Willmore_check}
\langle \text{Re}\,D_z\kappa,D_{\bar{z}}D_{\bar{z}}\kappa+\tfrac{1}{2}\bar{c}\kappa\rangle =
\text{Re}\,\langle D_z\kappa,D_{\bar{z}}D_{\bar{z}}\kappa+\tfrac{1}{2}\bar{c}\kappa\rangle \,,
\end{equation}
where we note that
$D_{\bar{z}}D_{\bar{z}}\kappa+\tfrac{1}{2}\bar{c}\kappa$ is real
by Codazzi's equation \eqref{eq:Codazzi}. Up to exact forms
\begin{align*}
-\langle D_z\kappa,D_{\bar{z}}D_{\bar{z}}\kappa+\tfrac{1}{2}\bar{c}\kappa\rangle \equiv
&\langle D_{\bar{z}}D_z\kappa,D_{\bar{z}}\kappa\rangle +\tfrac{1}{4}\bar{c}_z\langle \kappa,\kappa\rangle \\\equiv
&\langle R^{D}_{\bar{z}z}\kappa,D_{\bar{z}}\kappa\rangle +(\tfrac{3}{2}\langle D_{\bar{z}}\kappa,\bar{\kappa}\rangle 
+\tfrac{1}{2}\langle \kappa,D_{\bar{z}}\bar{\kappa}\rangle )\langle \kappa,\kappa\rangle \\\equiv
& \tfrac{3}{2}\langle R^{D}_{\bar{z}z}\kappa,D_{\bar{z}}\kappa\rangle \,.
\end{align*}
Here we used Gauss' equation \eqref{eq:Gauss}, the Ricci equation
\eqref{eq:Ricci} and the formula
\[
\langle \kappa,\kappa\rangle \langle \kappa,\bar{\kappa}\rangle _{\bar{z}}=-\langle R^{D}_{\bar{z}z}\kappa,D_{\bar{z}}\kappa\rangle 
+(3\langle D_{\bar{z}}\kappa,\bar{\kappa}\rangle 
+\langle \kappa,D_{\bar{z}}\bar{\kappa}\rangle )\langle \kappa,\kappa\rangle \,.
\]
Therefore we conclude that our third order flow, which we will call
the {\em Novikov--Veselov flow}, preserves the Willmore energy in case
flatness of the normal bundle is preserved. This is certainly the
case for surfaces in 3-space and also, as we will see later on, for
isothermic surfaces in any codimension.

At present we are unable to say more about the general codimension
case for the Novikov--Veselov flow. This is mainly due to our low
technology setup, which makes any discussion of the normal bundle
flow a notational debauch.

There is still much to do here: first we must see that our flow
coincides with that defined in Euclidean terms in the existing
literature \cite{KonTai95A,Tai97}.  Further a recursive scheme for
constructing higher flows is required.  For $n=3$, the first question
is answered in the affirmative by Richter \cite{Ric97} who proposes also
$\sigma_{2n+1}=\Re (D_z)^n\kappa$ for the higher flows.  We shall
return to these matters elsewhere.

\subsection{The Davey--Stewartson flows}\label{subsec:DS_flow}
In the previous section we indicated how to construct a hierarchy
of M\"obius invariant odd order flows on immersed tori in
$n$-space. To obtain even order flows one needs more structure in
the normal bundle, which is available for surfaces in $4$-space.
There, the normal bundle has a complex structure $J$ which we can
use to define a zero order flow, namely rotation in the normal
bundle by $J$. The corresponding deformation of the Hopf
differential is given by
\[
D_{t_0}\kappa=J\kappa
\]
and we take
\[
\sigma=\text{Re}\,J\kappa
\]
as our normal deformation for the construction of the second order
flow. To obtain the conformal deformation $f(t):M\to S^4$, we have
to solve for the tangential variation $b$ in \eqref{eq:Psi_dot}.
As before, this amounts to solving the $\delbar$-problem
$\bar{b}_z=2\langle \sigma,\kappa\rangle $. Now
\[
2\langle \sigma,\kappa\rangle =\langle J\kappa,\kappa\rangle +\langle J\bar{\kappa},\kappa\rangle =\langle J\bar{\kappa},\kappa\rangle 
\]
since $J$ is skew. But Ricci's equation \eqref{eq:Ricci} implies
\[
\int\langle J\bar{\kappa},\kappa\rangle =\pi\deg N_f M\,,
\]
so that the $\delbar$-problem \eqref{eq:conf_def_eq} is solvable
if and only if the normal bundle of $f$ has degree zero. This
condition is preserved under any continuous deformation, and we
obtain a unique second order flow by normalizing $\int b=0$.

This flow also preserves the Willmore energy \eqref{eq:W_dot}: up
to exact forms
\[
\langle \text{Re}\,J\kappa,D_{\bar{z}}D_{\bar{z}}\kappa+\tfrac{1}{2}\bar{c}\kappa\rangle =
\text{Re}\langle J\kappa,D_{\bar{z}}D_{\bar{z}}\kappa\rangle \equiv
\text{Re}\langle D_{\bar{z}}\kappa,JD_{\bar{z}}\kappa\rangle =0\,,
\]
where we used that $DJ=0$ and $J$ is skew.

Again there is more to be done, both in making contact with the
existing literature \cite{Kon00,KonLan00} and constructing higher
flows.

\subsection{Flows on isothermic surfaces}\label{subsec:isothermic}
As an explicit example of the above discussions, we apply the
Novikov--Veselov and Davey--Stewartson flow to isothermic surfaces
in $S^n$ \cite{Bur00}. Recall that $f:M\to S^n$ is isothermic if the Hopf
differential is real, i.e., if $\kappa=\bar{\kappa}$. Note that in
this case Ricci's equation \eqref{eq:Ricci} immediately implies
the flatness of the normal bundle. Thus we can define a ``vector"
version of the Novikov--Veselov flow and the Davey--Stewartson flow
in $4$-space on isothermic surfaces, both preserving the Willmore
energy, once we have shown that the isothermic condition is preserved under
the flows. 

First we note that in the isothermic case the $\delbar$-problem
can be explicitly solved for both flows. Since $\kappa$ is real,
the integrability equations for an isothermic surface become
\begin{gather*}
c_{\bar{z}}=4\langle \kappa,\kappa\rangle _z\,,\\
\Im\,(D_{\bar{z}}D_{\bar{z}}\kappa+\tfrac{1}{2}\bar{c}\kappa)=0\,,\\
R^{D}=0\,.
\end{gather*}
As discussed in section~\ref{subsec:NV_flow} the Novikov--Veselov
flow is (up to the factor of $2$) given by the normal variation
$\sigma=2\text{Re}\,D_z\kappa$. Therefore
\[
\bar{b}_z=2\langle \sigma,\kappa\rangle =\langle \kappa,\kappa\rangle _z+\langle \kappa,\kappa\rangle _{\bar{z}}
=(\langle \kappa,\kappa\rangle +\tfrac{\bar{c}}{4})_z
\]
and we obtain the explicit solution
\[
b=\langle \kappa,\kappa\rangle +\tfrac{c}{4}
\]
of the $\delbar$-problem \eqref{eq:conf_def_eq}. To see that the
Novikov--Veselov flow preserves isothermic surfaces, we need to
check whether $\kappa=\bar{\kappa}$ is preserved under the flow.
Recalling \eqref{eq:kappa_dot} and the fact that $\kappa$ is real,
we immediately calculate
\begin{align*}
\text{Im}(D_{t_3}\kappa)&=\text{Im}(D_z D_z D_z\kappa+D_z D_z
D_{\bar{z}}\kappa+\tfrac{c}{2}(D_z+D_{\bar{z}})\kappa+2b_z\kappa)\\
&=(D_z+D_{\bar{z}})\text{Im}(D_z
D_z\kappa+\tfrac{c}{2}\kappa)-\tfrac{1}{2}\text{Im}(
c_{\bar{z}}-4\langle \kappa,\kappa\rangle _z)\kappa\\
&=0\,,
\end{align*}
where we used the integrability equations. Therefore isothermic
surfaces are preserved under the Novikov--Veselov flow and
\[
D_{t_3}\kappa=\text{Re}(2D_z D_z
D_z\kappa+\tfrac{3}{2}cD_z\kappa+\tfrac{3}{4}c_z\kappa+
2(\langle \kappa,\kappa\rangle D_z\kappa-\langle \kappa,D_z\kappa\rangle \kappa))\,,
\]
which is a modified vector version of the usual Novikov--Veselov
equation found in the literature \cite{Tai97}. Notice that the last
term
\[
\langle \kappa,\kappa\rangle D_z\kappa-\langle \kappa,D_z\kappa\rangle \kappa=0
\]
for surfaces in $3$-space.

A similar, but simpler, computation shows that the
Davey--Stewartson flow preserves isothermic surfaces in $4$-space.
In this case the normal deformation
\[
\sigma=\text{Re}\,J\kappa=J\kappa
\]
since $\kappa$ is real. The $\delbar$-problem
\[
\bar{b}_z=2\langle \sigma,\kappa\rangle =2\langle J\kappa,\kappa\rangle =0
\]
is solved by $b=0$. Using again \eqref{eq:kappa_dot} we calculate
\[
\text{Im}D_{t_2}\kappa=\text{Im}(D_z D_z
J\kappa+\tfrac{c}{2}J\kappa)=0\,.
\]
This gives
\[
D_{t_2}\kappa=J(D_z D_z\kappa+\tfrac{c}{2}\kappa)
\]
for the Davey--Stewartson flow on $\kappa$ for isothermic surfaces
in $4$-space.

\providecommand{\bysame}{\leavevmode\hbox to3em{\hrulefill}\thinspace}
\providecommand{\MR}{\relax\ifhmode\unskip\space\fi MR }
\providecommand{\MRhref}[2]{%
  \href{http://www.ams.org/mathscinet-getitem?mr=#1}{#2}
}
\providecommand{\href}[2]{#2}

\end{document}